# A Generalization and Extension of an Autoregressive Model

Satheesh S

NEELOLPALAM, S. N. Park Road
Trichur – 680 004, **India.**
ssatheesh@sancharnet.in

Sandhya E

Department of Statistics, Prajyoti Niketan College
Pudukkad, Trichur – 680 301, **India.**
esandhya@hotmail.com

Rajasekharan K E

Department of Statistics, University of Calicut
Calicut University PO – 673 635, **India.**

**Abstract.** Generalizations and extensions of a first order autoregressive model of Lawrance and Lewis (1981) are considered and characterized here.

**1. Introduction.**

In this paper we consider the first order autoregressive (AR(1)) model (2.2) of Lawrance and Lewis (1981), generalize it and extend it to the maximum scheme giving stationary solutions to them. Here a sequence of *r.v*s $\{X_n, n>0 \text{ integer}\}$ defines the AR(1) scheme if for some $0<p<1$ there exists an innovation sequence $\{\varepsilon_n\}$ of *i.i.d r.v*s such that ;

$$
\begin{aligned}
X_n &= \varepsilon_n, & &\text{with probability } p \\
&= X_{n-1} + \varepsilon_n, & &\text{with probability } (1-p).
\end{aligned}
\qquad (1)
$$

Jose and Pillai (1995) showed that (1) is stationary and defined for each $0<p<1$ *iff* $X_n$ is geometrically infinitely divisible (ID). Recently Seetha Lekshmi and Jose (2006) has also discussed the model (1) and characterized the geometric Pakes law.

*Completed on 20 October 2006.*

This investigation is motivated by the possibility whereby the AR(1) sequence $\{X_n\}$ is composed of $k$ independent AR(1) sequences $\{Y_{i,n}\}$, $i = 1,2, \ldots, k$, and where for each $n>0$ integer $Y_{i,n}$, $i = 1,2, \ldots, k$ are identically distributed. That is, for each $n$, $X_n \stackrel{d}{=} \sum_{i=1}^{k} Y_{i,n}$ and $Y_{i,n}$, $i = 1,2, \ldots, k$ are i.i.d, $k$ being a fixed positive integer. For example, the variable $X_n$ could be the quantity of water collected in a dam in a week or the income from the sales of a particular item by an agency having $k$ different outlets. In these cases the resultant observation $X_n$ is either the sum of the quantities $Y_{i,n}$, $i = 1,2, \ldots, k$ (=7) of water collected in each of the 7 days or the sum of the incomes $Y_{i,n}$, $i = 1,2, \ldots, k$ from the sales in $k$ different outlets. Satheesh, *et al*. (2006) has generalized the AR(1) scheme of Gaver and Lewis (1980) on the above lines to characterize gamma-semi-stable and gamma-max-semi-stable laws. In this generalization of the AR(1) model we make a fruitful application of two approaches to the distribution of sums and maximums in samples of random size.

Here in section.2 we generalize the AR(1) model (1) where the observations $X_n$ are the sum of $k$ i.i.d r.vs, and then extend it to the maximum scheme in section.3. In section.4 we discuss another possible approach to the generalization.

**2. A Generalization of the AR(1) Additive Scheme.**

Harris($a,k$) law on $\{1, 1+k, 1+2k, \ldots\}$ is a generalization of the geometric law (when $k=1$) and is described by its probability generating function (PGF)

$$P(s) = \frac{s}{\{a - (a-1)s^k\}^{1/k}}, \quad k>0 \text{ integer and } a>1.$$

The stability properties of the Harris($a,k$) law in the minimum and maximum schemes were studied by Satheesh and Nair (2004). Distributional and divisibility properties, simulation and estimation of Harris($a,k$) law are discussed in Sandhya, *et al*. (2005).

Let us generalize (1) as follows. A sequence $\{Y_{i,n}\}$ defines a generalized AR(1) scheme if for some $0<p<1$ there exists an innovation sequence $\{\varepsilon_{i,n}\}$ of *i.i.d r.v*s such that;

$$\left.\begin{array}{l} \sum_{i=1}^{k} Y_{i,n} = \sum_{i=1}^{k} \varepsilon_{i,n}, \quad \text{with probability } p \\ \phantom{\sum_{i=1}^{k} Y_{i,n}} = \sum_{i=1}^{k} Y_{i,n-1} + \sum_{i=1}^{k} \varepsilon_{i,n}, \quad \text{with probability } (1-p) \end{array}\right\} \quad (2)$$

for all $n>0$ integer. An assumption here (and in similar structures to follow) is that the Bernoulli *r.v* deciding the exclusion of $\sum_{i=1}^{k} Y_{i,n-1}$ with probability $p$ is independent of $\{\varepsilon_{i,n}\}$.

In terms of characteristic functions (CF) and assuming stationarity (2) is equivalent to;

$$f_y^k(t) = p f_\varepsilon^k(t) + (1-p) f_y^k(t) f_\varepsilon^k(t)$$

$$= \frac{p f_\varepsilon^k(t)}{1-(1-p) f_\varepsilon^k(t)}. \text{ That is;}$$

$$f_y(t) = \frac{f_\varepsilon(t)}{\{a-(a-1)f_\varepsilon^k(t)\}^{1/k}}, \, a=\tfrac{1}{p}, \text{ Hence;} \quad (2a)$$

**Theorem.2.1** A sequence $\{Y_{i,n}\}$ defines the model (2) that is stationary for some $0<p<1$ *iff* $Y_{i,n}$ is a Harris$(a,k)$-sum and the innovation sequence is the components in the Harris$(a,k)$-sum, $a=\tfrac{1}{p}$.

Suppose we demand (2) to be satisfied for each $0<p<1$. That is, for each $0<p<1$ there exists a sequence $\{\varepsilon_{i,n,p}\}$ of *i.i.d r.v*s, so that (2a) is true for each $0<p<1$ with $f_\varepsilon(t)$ replaced by $f_{\varepsilon,p}(t)$. To discuss this we need the notion of *N*-ID laws of Gnedenko and Korolev (1996, *p*.137-152) that generalizes the geometrically ID laws as follows. Let $N = \{N_\theta, \theta \in \Theta\}$ be a family of positive integer valued *r.v*s having finite mean with PGF $\{P_\theta, \theta \in \Theta\}$, $\{\theta \in \Theta\}$ being the parameter space. A CF $f(t)$ is *N*-ID if for each $\theta \in \Theta$ there exists a CF $f_\theta(t)$ such that

$$f(t) = P_\theta\{f_\theta(t)\}, \, \forall \, t \in \mathbf{R}.$$

Of course, an underlying assumption is that the distributions corresponding to $P_\theta$ and $f_\theta$ are independent for each $\theta \in \Theta$. Their approach also required that the family of PGFs $\{P_\theta: \theta \in \Theta\}$ formed a commutative semi-group with respect to the operation of convolution. Then they showed that $f(t)$ is N-ID iff

$$f(t) = \varphi\{-\log h(t)\}$$

where $\varphi$ is a Laplace transform (LT) and $h(t)$ is the CF of an ID law. They also showed that $\{P_\theta: \theta \in \Theta\}$ formed a commutative semi-group iff $P_\theta$ and $\varphi$ are related by

$$P_\theta(z) = \varphi\{\varphi^{-1}(s)/\theta\}, \theta \in \Theta.$$

Now, if (2a) is to be true for each $0<p<1$ with $f_\varepsilon(t)$ replaced by $f_{\varepsilon,p}(t)$, what we require is that $f_y(t)$ is Harris-ID. Note that if $\varphi(s) = \dfrac{1}{(1+s)^{1/k}}$, the LT of the gamma law, then, $\varphi^{-1}(s) = \dfrac{1-s^k}{s^k}$ and $\varphi\{\varphi^{-1}(s)/p\} = \dfrac{s}{\{a-(a-1)s^k\}^{1/k}}$, $a = \dfrac{1}{p}$ is the PGF of the Harris($a,k$) law. Hence;

**Theorem.2.2** A CF $f(t)$ is Harris-ID iff $f(t) = \dfrac{1}{(1-\log h(t))^{1/k}}$, $k>0$ integer, where $h(t)$ is some CF that is ID.

Though a corollary to theorem.2.2, in the context of the model (2) we have;

**Theorem.2.3** $\{Y_{i,n}\}$ defines the model (2) that is stationary for each $0<p<1$, iff $Y_{i,n}$ is Harris($a,k$)-ID, $a = \dfrac{1}{p}$.

**3. The AR(1) Maximum Scheme and its Generalization.**

The max-analogue of the AR(1) model of Gaver and Lewis (1980) was introduced and characterized in Satheesh and Sandhya (2006). The max-analogue of (1) is

$$X_n = \varepsilon_n, \quad \text{with probability } p$$
$$= X_{n-1} \vee \varepsilon_n, \quad \text{with probability } (1-p). \tag{3}$$

In terms of *d.fs* and assuming stationarity (3) is equivalent to;

$$F(x) = pF_\varepsilon(x) + (1-p) F(x) F_\varepsilon(x)$$
$$= \frac{pF_\varepsilon(x)}{1-(1-p)F_\varepsilon(x)}, \text{ which proves;} \tag{3a}$$

**Theorem.3.1** $\{X_n\}$ defines the model (3) that is stationary for some $0<p<1$ *iff* the *d.f* of $X_n$ is a geometric maximum and the innovation sequence has the distribution of the components.

Suppose we demand (3) to be satisfied for each $0<p<1$, then (3a) is true for each $0<p<1$ with $F_\varepsilon(x)$ replaced by $F_{\varepsilon,p}(x)$. Hence by Rachev and Resnick (1991) $X_n$ must be max-geometrically-ID. Hence we have;

**Theorem.3.2** $\{X_n\}$ defines the model (3) that is stationary for each $0<p<1$ *iff* $X_n$ is max-geometrically-ID with d.f $F(x) = \dfrac{1}{1-\log H(x)}$, where $H(x)$ is a *d.f*.

**Remark.3.1** Strictly speaking, $H(x)$ above should be max-ID, but then in the univariate case every *d.f* is max-ID.

Next consider the generalization of (3) as done in (2) as follows.

$$\bigvee_{i=1}^{k} Y_{i,n} = \bigvee_{i=1}^{k} \varepsilon_{i,n} \quad \text{with probability } p$$
$$\bigvee_{i=1}^{k} Y_{i,n} = \{\bigvee_{i=1}^{k} Y_{i,n-1}\} \vee \{\bigvee_{i=1}^{k} \varepsilon_{i,n}\} \text{ with probability } (1-p) \tag{4}$$

In terms of *d.fs* and assuming stationarity we have (as done to arrive at (2a));

$$F(x) = \frac{F_\varepsilon(x)}{\{a-(a-1)F_\varepsilon^k(x)\}^{1/k}}, \ a = \tfrac{1}{p}. \text{ This proves;} \tag{4a}$$

**Theorem.3.3** A sequence $\{Y_{i,n}\}$ defines the model (4) that is stationary for some $0<p<1$ *iff* the *d.f* of $Y_{i,n}$ is a Harris($a,k$)-maximum and the innovation sequence has the distribution of the components in the Harris($a,k$)-maximum, $a=\frac{1}{p}$.

Suppose we demand (4) to be satisfied for each $0<p<1$, then (4a) is true for each $0<p<1$ with $F_\varepsilon(x)$ replaced by $F_{\varepsilon,p}(x)$. To discuss this we need the max-analogue of the notion of *N*-ID laws in Gnedenko and Korolev (1996) which we develop now.

**Definition.3.1** Let $\varphi$ be a LT and $N_\theta$ a positive integer-valued *r.v* having finite mean with PGF $P_\theta(s) = \varphi(\frac{1}{\theta}\varphi^{-1}(s))$, $\theta \in \Theta$. A *d.f* $F$ is *N*-max-ID if for each $\theta \in \Theta$ there exists a *d.f* $F_\theta$ that is independent of $N_\theta$ such that $F(x) = P_\theta(F_\theta(x))$, $\forall x \in \mathbf{R}$.

**Theorem.3.4** The limit of a sequence of *N*-max-ID laws is *N*-max-ID.

*Proof.* Follows as in the proof of property.4.6.2 in Gnedenko and Korolev (1996, *p*.145).

**Theorem.3.5** Let $\varphi$ be a LT, $a_n>0$ are constants and $G_n$ are *d.f*s. A *d.f* $F$ is *N*-max-ID *iff*

$$F(x) = \operatorname*{Lt}_{n \to \infty} \varphi(a_n(1-G_n(x))), \forall x \in \mathbf{R}. \tag{5}$$

*Proof.* We closely follow the proof of theorem.4.6.2 in Gnedenko and Korolev (1996, *p*.146). The "*if* part" follows by noting that for each $a>0$, $F_a(x) = \varphi(a(1-G(x)))$ is a *d.f* being a mixture of Poisson maximums of $G(x)$'s where the mixing law has LT $\varphi$. Further,

$$F_a(\mathrm{x}) = \varphi\{\tfrac{1}{\theta}\varphi^{-1}[\varphi(a\theta(1-G(x)))]\} = \varphi\{\tfrac{1}{\theta}\varphi^{-1}[F_{a\theta}(x)]\} \text{ for each } \theta \in \Theta.$$

Thus $F_a(\mathrm{x})$ is *N*-max-ID or the *d.f* under the limit sign in (5) is *N*-max-ID. Hence by theorem.3.4 the limit *d.f* is again *N*-max-ID.

To prove the converse, let the *d.f* $F(x)$ be an *N*-max-ID. Then by definition.3.1 it can be represented as $F(x) = P_\theta(F_\theta(x)) = \varphi(\tfrac{1}{\theta}\varphi^{-1}(F_\theta(x)))$, for each $\theta \in \Theta$ and $\forall x \in \mathbf{R}$. Hence $F_\theta(x) = \varphi(\theta\varphi^{-1}(F(x)))$ is also a *d.f* for each $\theta \in \Theta$.

Setting $a = \frac{1}{\theta}$ and $G(x) = \varphi(\theta\varphi^{-1}(F(x)))$ in the arguments used in the proof of the "*if* part" we can see that the *d.f*

$$\varphi(\tfrac{1}{\theta}(1-G(x))) = \varphi\{\tfrac{1}{\theta}[1-\varphi(\theta\varphi^{-1}(F(x)))]\}$$

is N-max-ID and by the properties of the function $\varphi$ (a LT)

$$F(x) = \underset{\theta \to 0}{Lt}\ \varphi\{\tfrac{1}{\theta}[1-\varphi(\theta\varphi^{-1}(F(x)))]\}\ \text{which completes the proof.}$$

**Theorem.3.6** A *d.f* $F$ is *N*-max-ID iff $F(x) = \varphi\{-\log H(x)\}$, where $\varphi$ is a LT and $H$ a *d.f*.

*Proof.* Follows by noting the remark.3.1 and theorem.3.5. Thus we have:

**Theorem.3.7** A *d.f* $F$ is Harris-max-ID iff $F(x) = \dfrac{1}{\{1-\log H(x)\}^{1/k}}$, $k>0$ integer, and $H$ a *d.f*. Proof follows as that of theorem.2.2.

**Theorem.3.8** $\{Y_{i,n}\}$ defines the model (4) that is stationary for each $0<p<1$ iff $Y_{i,n}$ is Harris($a,k$)-max-ID with *d.f* $F(x) = \dfrac{1}{\{1-\log H(x)\}^{1/k}}$, for a *d.f* $H(x)$, $a=\tfrac{1}{p}$.

**4. Another Look at the Models (2) and (4).**

A problem in the above development is how to identify the CF $h$ or the *d.f* $H$. This motivates us to take another look at the models (2) and (4) by using the notions of φ-ID and φ-max-ID laws from Satheesh (2002, 2004). In the terminology used in the above works geometrically-ID and Harris-ID laws are respectively exponential-ID and gamma-ID laws. Similarly their max-analogues are exponential-max-ID and gamma-max-ID laws. This is because the terminology is derived from the form of the limiting LT of the *r.v* $\theta N_\theta$ as $\theta \downarrow 0$. The discussions in the above works are based on the following two lemmas and definition and invoking transfer theorems for random ($N_\theta$) sums and maximums.

**Lemma.4.1** $\wp_\varphi = \{P_\theta(s) = s^j \varphi\{(1-s^k)/\theta\},\ 0<s<1,\ j\geq 0\ \&\ k\geq 1\ \text{integer and}\ \theta>0\}$ describes a class of PGFs for any given LT $\varphi$.

**Lemma.4.2** Given a r.v $U$ with LT $\varphi$, the integer valued r.vs $N_\theta$ with PGF $P_\theta$ in the class $\wp_\varphi$ described in lemma.4.1 satisfy

$$\theta N_\theta \xrightarrow{d} kU \text{ as } \theta \downarrow 0.$$

**Definition.4.1** Let $\varphi$ be a LT. A CF $f$ (d.f $F$) is $\varphi$-ID ($\varphi$-max-ID) if for each $\theta \in \Theta$, there exists a CF $h_\theta$ (d.f $H_\theta$), a PGF $P_\theta \in \wp_\varphi$ that is independent of $h_\theta$ ($H_\theta$), such that $P_\theta\{h_\theta\} \to f$ ($P_\theta\{H_\theta\} \to F$) as $\theta \downarrow 0$ through a $\{\theta_n\}$.

Assuming the existence of the classical limit of sums (maximums) of i.i.d r.vs $\{X_{\theta,j}\}$ as $\theta \downarrow 0$ we have the following result invoking the transfer theorems for sums (maximums).

The limit law of $N_\theta$-sums ($N_\theta$-maximums) of i.i.d r.vs $\{X_{\theta,j}\}$ as $\theta \downarrow 0$, where the PGF of $N_\theta$ is a member of $\wp_\varphi$, is necessarily $\varphi$-ID ($\varphi$-MID). Conversely, for any LT $\varphi$ the $\varphi$-ID ($\varphi$-MID) law can be obtained as the limit law of $N_\theta$-sums ($N_\theta$-maximums) of i.i.d r.vs as $\theta \downarrow 0$ for each $N_\theta$ whose PGF is in $\wp_\varphi$. The $N_\theta$-sum version here is theorem.2.8 in Satheesh (2004).

Now, instead of demanding (4) to be satisfied for each $0<p<1$, let us consider the apparently weaker requirement that;

$$F(x) = \underset{p\downarrow 0}{Lt} \frac{pF_{\varepsilon,p}(x)}{\{1-(1-p)F_{\varepsilon,p}^k(x)\}^{1/k}}, \text{ where } p\downarrow 0 \text{ through a sequence } \{p_n\}. \tag{6}$$

Notice that the Harris PGF is $P(s) = s\varphi(\frac{1}{\theta}(1-s^k))$, where $\varphi(s)$ is the LT $\frac{1}{(1+s)^{1/k}}$, $\theta = \frac{p}{1-p}$, and $p=\frac{1}{a}$. Hence by the two lemmas above and the transfer theorem for maximums, the Harris maximums in (6) converge weakly to the d.f $\frac{1}{\{1+k(-\log H(x))\}^{1/k}}$, provided $H(x)$ is the weak limit of $F_{\varepsilon,p}$'s as $p \downarrow 0$ through $\{p_n\}$. $H(x)$ being univariate we only need that $H(x)$ is a proper d.f so that it is max-ID. Hence we have:

**Theorem.4.1** $\{Y_{i,n}\}$ defines the model (4) that is stationary for a null sequence $\{p_n\}$ iff $Y_{i,n}$ is gamma-max-ID with d.f $\dfrac{1}{\{1+k(-\log H(x))\}^{1/k}}$, where $H(x)$ is the weak limit of $\{F_{\varepsilon,p}\}$ as $p\downarrow 0$ through $\{p_n\}$.

Similarly one has the following result for the model (2).

**Theorem.4.2** $\{Y_{i,n}\}$ defines the model (2) that is stationary for a null sequence $\{p_n\}$ iff $Y_{i,n}$ is gamma-ID with CF $\dfrac{1}{\{1+k(-\log h(t))\}^{1/k}}$, where $h(t)$ is the limiting CF of $\{f_{\varepsilon,p}\}$ as $p\downarrow 0$ through $\{p_n\}$.

**Remark.4.1** Whether we approach the problem, as done in sections 2 and 3 above or as in this section we need consider weak limits of distributions concerned. The advantage of using the approach in this section is that the requirement is apparently weaker and the limit law is a function of the CF $h$ (or d.f $H$) which is the weak limit of the distributions of the innovation sequence described by their CFs $\{f_{\varepsilon,p}\}$ in (2) (or their d.fs $\{F_{\varepsilon,p}\}$ in (4)). This advantage will not be there if we approach the problem as in Gnedenko and Korelev (1996) and its max-analogue. Notice that gamma-ID and gamma-max-ID laws identified in theorems 4.2 and 4.1 respectively, also satisfy the requirement of Harris($a,k$)-sum and Harris($a,k$)-maximum for each $0<p<1$, $p=\frac{1}{a}$.

**References.**

**Gaver, D P and Lewis, P A W (1980),** First-order auto regressive gamma sequences and point processes, *Adv. Appl. Prob.* **12**, 727–745.

**Gnedenko, B V (1982),** On limit theorems for a random number of random variables, *Lecture Notes in Mathematics*, **Vol.1021**, Fourth USSR-Japan Symp. Proc., Springer, Berlin, 167-176.

**Gnedenko, B V and Korolev, V Yu. (1996),** *Random Summation - Limit Theorems and Applications*, CRC Press, Boca Raton.

**Jose, K K and Pillai, R N (1995),** Geometric infinite divisibility and its applications in autoregressive time series modeling, in *Stochastic Process and Its Applications*, Editor V Thankaraj, Wiley Eastern, New Delhi.

**Lawrance, A J and Lewis, P A W (1981),** A new autoregressive time series modeling in exponential variables (NEAR(1)), *Adv. Appl. Probab.* **13**, 826–845.

**Rachev, S T and Resnick, S (1991),** Max geometric infinite divisibility and stability, *Comm. Statist. - Stoch. Models* **7**, 191–218.

**Seetha Lekshmi, V and Jose, K K (2006),** Autoregressive processes with Pakes and geometric Pakes generalized Linnik marginals, *Statist. Probab. Lett.* **76**, 318-326.

**Sandhya, E; Sherly, S and Raju, N (2005),** Harris family of discrete distributions, accepted to appear in the *J. Kerala Statist. Assoc.*, 2007.

**Satheesh, S (2002),** Aspects of randomization in infinitely divisible and max-infinitely divisible laws, *ProbStat Models* **1**, 7–16.

**Satheesh, S (2004),** Another look at random infinite divisibility. *Statist. Meth.* **6**, 123-144.

**Satheesh, S and Nair, N U (2004),** On the stability of geometric extremes, *J. Ind. Statist. Assoc.* **42**, 99–109.

**Satheesh, S and Sandhya, E (2006),** Max-semi-selfdecomposable laws and related processes, *Statist. Probab. Lett.* **76**, 1435–1440.

**Satheesh, S; Sandhya, E and Sherly, S (2006),** A generalization of stationary AR(1) schemes, accepted to appear in *Statist. Meth.*, 2007.